# Identifying Long Term Voltage Stability Caused by Distribution Systems vs Transmission Systems

Amarsagar Reddy Ramapuram M.
Department of Electrical and Computer Engineering
Iowa State University, Ames, IA,
amar@iastate.edu

Ankit Singhal
Department of Electrical and Computer Engineering
Iowa State University, Ames, IA,
ankit@iastate.edu

Venkataramana Ajjarapu
Department of Electrical and Computer Engineering
Iowa State University, Ames, IA,
vajjarapu@iastate.edu

*Abstract* — **Monitoring the long term voltage stability of the power grid is necessary to ensure its secure operation. This paper presents a new phasor based methodology that distinguishes between long term voltage stability caused by distribution systems versus transmission systems. From a conceptual understanding of a simplified system, a Transmission-Distribution Distinguishing Index (TDDI) is proposed to distinguish between the two scenarios. A methodology to calculate the TDDI for multi-bus systems using quasi-steady state phasor measurements is described and validating results are presented for the IEEE 9 Bus system with a load replaced by various distribution feeders. The results verify that the TDDI can indeed be used to distinguish between transmission limited and distribution limited systems. This information can be utilized by the operator to effectively choose controls in distribution and transmission systems to improve the system margin.**

## I. Introduction

There is increasing pressure on power system operators and on electric utilities to utilize the existing grid infrastructure to the maximum extent possible and this mode of operation can lead to long term voltage stability problems. To handle this, operators are adopting real-time tools using Wide-Area measurements (WAMS) and Phasor Measurement Units (PMUs) that are providing them with better situational awareness. The increasing number of PMUs in the grid have led to various online Voltage Stability Indices (VSI) being proposed in recent times to monitor the grid in real time [1].

Traditionally, the VSI's were calculated at a bus by estimating the Thevenin Equivalent using local PMU Voltage and Current measurements at a Bus [2-6]. However, all these methods assume an aggregated load at the transmission level and do not consider the sub-transmission system or the distribution system where the loads are actually present. Ignoring the distribution feeder network and the distribution of loads in the sub-transmission/ distribution network will lead to an error in the voltage stability assessment.

Furthermore, as the distribution systems are often operated close to their limits (for economic reasons), considering their topology and loads into the voltage stability assessment might provide insights to operators and planners on how to improve the system behavior. In fact, voltage collapse in distribution feeders has been identified as a critical issue for some time [7] and a major blackout in 1997 in the S/SE Brazilian system is attributed to a voltage instability problem in one of the distribution feeders that spread to the transmission grid [8].

Recently, techniques incorporating the distribution system in the transmission system analysis have been proposed [9] and have been utilized to verify how the increase in Distributed Generation (DG) can improve both the overall system margin [10] and the distribution system margin [11]. However, as far as the authors know, none of the existing methods distinguish between the voltage stability caused by distribution system versus transmission systems. Our previous paper [12] describes a method to determine if the voltage stability limit is due to the distribution system or the transmission system. This method performs a continuation power flow [13] and compares the resultant nose point with a predetermined hypersurface based on the distribution topology. The methodology requires a full-fledged CPF routine along with the calculation of the hyper surfaces, making it time consuming for power system operations and so an online methodology would be preferred.

In this paper we address this issue by presenting a technique based on phasor measurements to estimate the voltage stability and to also determine if the limit is due to the transmission or distribution systems. This information will be useful for operators, especially as the control of DG devices in distribution limited systems can lead to a larger percentage increase in the margin [12] and a load shedding action on the distribution limited systems will lead to a larger improvement in system margin [12]. Thus, determination of the limiting system can be used to improve voltage stability with minimum control. This paper starts by describing a conceptual understanding of the methodology on a simple system (Section II), presents a technique to estimate the parameters of the equivalent circuit for multi-bus systems (Section III), describes how the method works on a test-system and compares it to existing results (Section IV) and finally concludes in Section V.

## II. Conceptual Understanding of method

A block diagram of the conventional power system is shown in Fig. 1, with the various generation, transmission and distribution circuits. The loads are in the distribution feeders and vary based on time of day, etc.

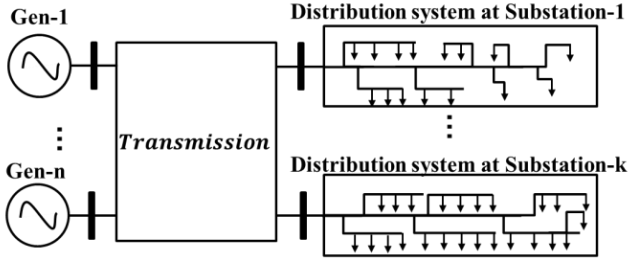

Fig. 1. A conventional power system topology showing interaction between the generators, the transmission system and the various distribution systems.

If the load in the system in Fig. 1 is increased at all buses along a particular load increase direction, then there is a load limit after which there is no solution possible [13]. This is the critical point of the system with respect to the long term voltage stability and occurs due to the limitations of the underlying transmission and distribution network. However, it is not straightforward to estimate which system (transmission or distribution) is limiting the critical load.

In order to distinguish between the systems limited by transmission and distribution network, a conceptual understanding of the phenomenon is necessary. The simplest system that is possible with a transmission and distribution system is shown in Fig 2. It is an extension of the standard Thevenin equivalent with an extra impedance to represent the distribution network. $E_{th}$ is the Thevenin voltage, $Z_T$ is the transmission system contribution to the Thevenin impedance, $Z_D$ is the distribution system contribution to the Thevenin impedance and $Z_L$ is the load impedance.

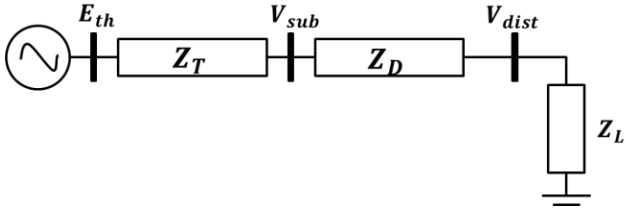

Fig. 2. A very simple transmission and distribution system.

For this particular system, the transmission and distribution are very simple and a straightforward method would be to compare the magnitude of the impedances of the transmission and distribution networks. If the transmission impedance is more than the distribution impedance ($|Z_T| > |Z_D|$), then the transmission network is the limit. Instead, if the distribution impedance is more than the transmission impedance ($|Z_D| > |Z_T|$), then the distribution network is the limit. The intuitive explanation is that the network which has a larger voltage drop is the main limiting network for voltage stability. Another way to look at it is as follows: if reducing the impedance of the transmission networks leads to a better voltage improvement than by reducing the impedance of the distribution networks by the same proportion, then the transmission system is the limiting factor.

Hence, for this simplified system, the ratio between the impedances can be used as a way to distinguish between transmission and distribution limited systems. However, instead of directly using the ratio, the following Transmission-Distribution Distinguishing Index (TDDI) is defined in Eq. (1)

$$TDDI = \log\left(\frac{|Z_T|}{|Z_D|}\right) \quad (1)$$

If the ratio $|Z_T|/|Z_D|$ is greater than 1 (transmission limited), then the value of TDDI is positive; if the ratio is less than 1 (distribution limited), the value of TDDI is negative; and if the ratio is equal to 1, the value of TDDI is zero. The reason for using the logarithm function in Eq. (1) can be understood from the following example. If the ratio $|Z_T|/|Z_D|$ is 3, then the transmission system contributes three times more than the distribution system and if the ratio is 1/3, then the distribution contributes three times more than the transmission. However, these are unequally far from the value when both transmission and distribution contribute equally ($|Z_T|/|Z_D|$ =1). Thus, this is a skewed metric and this is resolved by the logarithm function. As $\log(x) = -\log(1/x)$, the TDDI of these two scenarios is equidistant from the case when both transmission and distribution contribute equally to the limit (TDDI=0). Now that the simple circuit has been analyzed conceptually, applying this method to a multi bus system requires a way to estimate the equivalent circuit parameters and this is presented in the next section.

### III. ESTIMATION OF PARAMETERS FROM MEASUREMENTS

As there is only a single load present in the equivalent circuit analyzed in the previous section, it implies that an equivalent circuit can be formed for every single load in the integrated transmission-distribution system. In order to estimate the parameters of the Thevenin equivalent circuit for a multi-bus system, conventional methods utilize phasor measurements of the load and the quasi static behavior of the system [2]. Utilizing these measurements, the total Thevenin impedance ($Z_T + Z_D$) can be estimated but not individual transmission equivalent ($Z_T$) or distribution equivalent ($Z_D$). An additional phasor measurement is necessary at the substation where the distribution feeder connects to the transmission system in order to estimate the values $Z_T$ and $Z_D$ separately.

Let $V_{sub}^{(i)}, V_{dist}^{(i)}$ be the voltage phasor measurement at the substation bus and distribution load bus at a time instant $i$. Let the current phasor measurement of the distribution load at instant $i$ be denoted by $I_{dist}^{(i)}$. As the equivalent in Fig. 2 is separately determined for every single load bus, the voltage and current phasor measurements are related to the equivalent circuit parameters through equations (2)-(4) at every time instant.

$$E_{th} = V_{sub}^{(i)} + I_{dist}^{(i)} \cdot Z_T \quad (2)$$

$$V_{sub}^{(i)} = V_{dist}^{(i)} + I_{dist}^{(i)} \cdot Z_D \quad (3)$$

$$Z_L = \frac{V_{dist}^{(i)}}{I_{dist}^{(i)}} \quad (4)$$

Equations (2) & (3) can be written in a matrix form as Eq. (5)

$$\begin{bmatrix} -1 & I_{dist}^{(i)} & 0 \\ 0 & 0 & I_{dist}^{(i)} \end{bmatrix} \begin{bmatrix} E_{Th} \\ Z_T \\ Z_D \end{bmatrix} = \begin{bmatrix} -V_{sub}^{(i)} \\ V_{sub}^{(i)} - V_{dist}^{(i)} \end{bmatrix} \quad (5)$$

This is a rank deficient set of equations and so measurements from at least 2 time instants are necessary for the estimation of the circuit parameters. Thus, a necessary assumption is that the circuit equivalent parameters will remain same for these time instants, which is a similar assumption made for conventional Thevenin methods [4]. The more the measurement instants, the larger the dimension of the matrix in Eq (5) and a least square estimate for the parameters $E_{th}, Z_T, Z_D$ can be performed. In practice the noise in the measurements can be handled by the robustness of the least square estimate. In this paper, the noise is ignored and so, the parameters can be estimated from phasor measurements at 2 instants with the expressions for $E_{th}, Z_T$ & $Z_T$ explicitly written as equations (6) - (8).

$$E_{th} = \left(\frac{V_{sub}^{(2)} \cdot I_{dist}^{(1)} - V_{sub}^{(1)} \cdot I_{dist}^{(2)}}{I_{dist}^{(1)} - I_{dist}^{(2)}}\right) \quad (6)$$

$$Z_T = -\frac{V_{sub}^{(1)} - V_{sub}^{(2)}}{I_{dist}^{(1)} - I_{dist}^{(2)}} \quad (7)$$

$$Z_D = \left(\frac{V_{sub}^{(1)} - V_{dist}^{(1)}}{I_{dist}^{(1)}} + \frac{V_{sub}^{(2)} - V_{dist}^{(2)}}{I_{dist}^{(2)}}\right) \cdot \frac{1}{2} \quad (8)$$

The load impedance can also be written as Eq. (9), which is the mean of the load impedance at the two instants.

$$Z_L = \left(\frac{V_{dist}^{(1)}}{I_{dist}^{(1)}} + \frac{V_{dist}^{(2)}}{I_{dist}^{(2)}}\right) \cdot \frac{1}{2} \quad (9)$$

Thus, using the phasor measurements at two instants, the equivalent circuit parameters can be estimated and utilized to distinguish between transmission limited and distribution limited systems. To measure the voltage and current phasors in the distribution network, MicroPMUs at low voltage circuits (12 kV, 33 kV, etc.) are necessary and at present they are being investigated by utilities [14], national labs [14] and universities [15] on how they improve system operation. In the future, we expect a few MicroPMUs to be deployed in key nodes in the distribution feeders and their data can be utilized for the proposed method.

As a different equivalent circuit is estimated for each load bus in the system, each load bus will have a corresponding TDDI. Thus, we first have to locate the load which is limiting voltage stability and then look at the TDDI for this particular bus. In order to detect the load which is limiting the instability, the conventional voltage stability index (VSI) is used [1]. The VSI is calculated from the circuit equivalent at every load bus and is given by Eq. (10)

$$VSI = \frac{|Z_T + Z_D|}{|Z_L|} \quad (10)$$

The closer the value of VSI to 1, the closer the load is to instability and so the load with the highest VSI is determined as the critical load. Once this load bus is determined, the equivalent circuit of this bus is used to calculate the TDDI and determine if the system is transmission limited or distribution limited. Results demonstrating this on multi-bus networks are described in the next section.

## IV. NUMERICAL VALIDATION OF METHODOLOGY

To test the methodology, an integrated Transmission-Distribution system is constructed from the IEEE 9 bus transmission system. The load at bus 5 (90 MW and 30 MVAR) is replaced by attaching multiple distribution feeder configurations in parallel. The distribution feeder configurations used have the topology shown in Fig. 3 and are based on the IEEE 4 bus distribution test system [16].

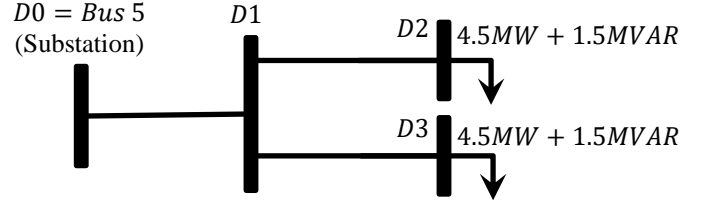

Fig. 3. The distribution topology used for validation. Bus D0 is the high voltage transmission bus (Bus 5 in this case). The Appendix has details on the line impedances and the loads at Bus D2 and D3.

The system is assumed to be balanced and the impedances of the various lines in the distribution system and the loads at Bus D2 and D3 are detailed in the Appendix for the various cases analyzed. As the system load of 90 MW is too large to be handled by any distribution system, multiple feeder configurations (specifically 10) are connected in parallel at bus 5 to ensure that the final system load is the same. The final integrated Transmission-Distribution system has a total of 39 buses (30 distribution buses and 9 transmission buses).

For this study, the power is increased in proportion to the original loading at all the loads (keeping power factor constant) and generators till the critical point is reached. MATPOWER[17] is used to run the CPF and analyze the results for the scenarios. Voltage and current phasors from consecutive points of the PV curve, which correspond to quasi-steady state measurements from the power system, are used to calculate the equivalent circuit for the loads in the distribution feeders. This equivalent is then used to calculate the VSI and the TDDI at every load bus. Since 10 identical feeder configurations are connected in parallel, analyzing the behavior of loads in a single feeder configuration is sufficient.

Two kinds of feeder configurations are used for this study to present the contrast between the transmission limited and distribution limited systems. Feeder configuration 1 (FC1) has feeders with large impedances and Feeder configuration 2 (FC2) has feeders with small impedances. This study has been previously reported in [12] where the voltage stability margin is calculated for three cases – (a) the load at bus 5 remains as specified by the standard IEEE 9 bus test system, (b) load at bus 5 replaced by 10 parallel FC1 and (c) load at bus 5 replaced by 10 parallel FC2. Table I lists the system voltage stability margin for these scenarios and Fig. 4 plots the PV curves for these scenarios.

TABLE I. TRANSMISSION SYSTEM LOAD

| Bus 5 Load | System Margin | Comment |
|---|---|---|
| Standard | 467.5 MW | - |
| FC1 | 163 MW | Large reduction in margin |
| FC2 | 419 MW | Small reduction in margin |

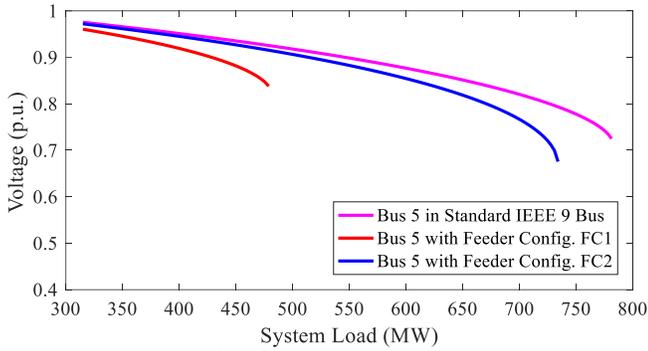

Fig. 4. PV curve at Bus 5 with 10 parallel connections of feeder configurations FC1 and FC2. Observe that the voltage at the critical point is ~0.71 p.u. for the standard case, ~0.85 p.u for FC1 and ~0.67 p.u. for FC2.

The large reduction in the margin due to replacing the load at bus 5 with FC1 is used as a metric along with hyper-surfaces in [12] to conclude that the distribution system is the cause for the long term voltage stability for this case. Similarly, [12] concludes that the small reduction in the margin with FC2 implies a limit in the transmission system. For the same systems, we will verify if the proposed method using TDDI also gives the same conclusions without needing to calculate the margin.

### A. Replacing Bus 5 load with Feeder Configuration FC1

The loads in feeder configuration FC1 are either at D2 or D3 in the distribution feeder and it is necessary to determine the critical load. To determine this, the VSI is calculated at D2 and D3 using Eq. (9) after replacing the load at bus 5 with FC1. Fig. 5 below plots the VSI versus the system load at buses D2 and D3. It can be observed that the VSI at bus D3 is higher than the VSI at bus D2 at all the load levels and more specifically at the critical load and so the critical bus is D3.

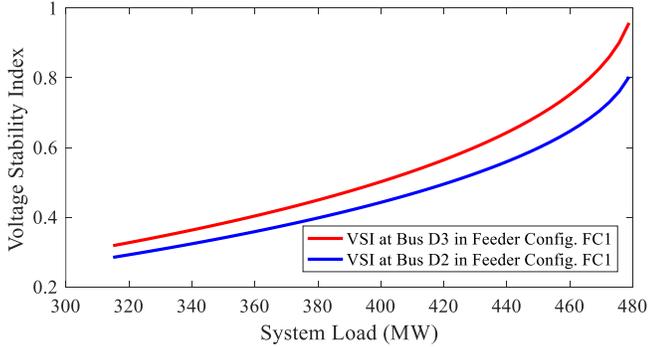

Fig. 5. VSI at D2 and D3 in feeder configuration FC1 versus system load. The VSI of D3 is greater than D2 and so it is the critical bus.

Now that the location of the critical load is determined, the TDDI at this bus is used to determine if the system is limited by distribution or by transmission networks. Fig. 6 plots the TDDI at Bus D2 and D3 versus the system load and it can be observed that the value of TDDI at bus D3 at critical loading is -0.4 which implies that $|Z_D| = 1.5 \cdot |Z_T|$ and the overall system is distribution limited. The CPF and Hyper-plane based method in [12] also has the same conclusion.

From Fig. 6, it can also be seen that the TDDI at bus D2 is close to 0, which seems to suggest that the transmission and distribution are equally limiting the load increase. So, if only this location is monitored, then a misleading conclusion can be made. Thus, it is important to calculate the TDDI at the critical bus when deciding if a system is transmission limited or distribution limited. Also, the TDDI at low loading (say 380 MW) is also negative with a value around -0.3. Thus, the TDDI at medium loading gives some information about the limiting system at the critical loading, but this might not always be the case.

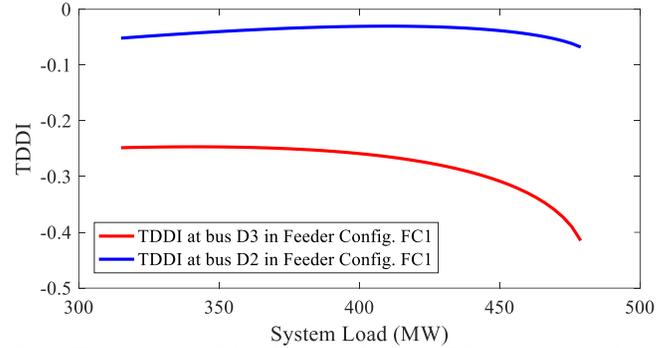

Fig. 6. TDDI at D2 and D3 in feeder configuration FC1 versus system load. The TDDI at bus D3 at critical loading is -0.4.

### B. Replacing Bus 5 load with Feeder Configuration FC2

Similar to feeder configuration FC1, the loads in FC2 are also at D2 or D3 in the distribution feeder and the VSI is calculated at D2 and D3 using Eq. (9) to determine the critical bus. Just as in the FC1, the critical bus in FC2 is also D3 and the VSI plot is omitted in the interest of space. The main difference between FC1 and FC2 is that the impedances are reduced and this reduction in the impedances improves the maximum loadability of the system to 734 MW.

To determine if the system is limited by distribution or by transmission networks, the TDDI is calculated and Fig. 7 plots the TDDI at Bus D2 and D3 versus the system load for FC2. It can be observed that the value of TDDI at bus D3 at critical loading is 0.71 which implies that $|Z_T| = 2 \cdot |Z_D|$ and the overall system is transmission limited which is the same conclusion in [12].

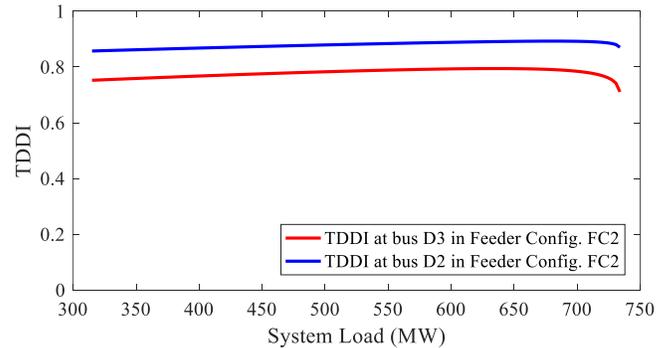

Fig. 7. TDDI at D2 and D3 in feeder configuration FC2 versus system load. The TDDI at bus D3 at critical loading is 0.71.

From Fig. 7, it can also be seen that the TDDI at bus D2 is close to 0.85, and so in this scenario the information about the limiting network can be approximately estimated from measuring the voltages at D2. Also, in this scenario, the TDDI is almost constant for a majority of the loading conditions and so the TDDI at the loading limit can be reasonably estimated from a medium loading condition.

From these results, we can say that the proposed method can indeed identify if the voltage stability limit is being caused by the distribution system or the transmission system directly from quasi-steady state phasor measurements.

## V. CONCLUSION AND FUTURE STUDIES

In this paper, a phasor measurement based methodology that distinguishes between distribution limited and transmission limited long term voltage stability is presented. An intuitive explanation, utilizing a simplified circuit is provided for the method and a technique to estimate the simplified equivalent circuit for multi-bus systems is described using phasor measurements from various time instants is described. Results for the IEEE 9 Bus system with a load being replaced by a set of feeders are presented. These results were validated with existing methodology based on continuation power flow [12] and the proposed method is able to distinguish between systems that are transmission limited and distribution limited. This property will enable the operators to quickly choose between various controls (e.g. DG in a particular distribution system, changing taps in the transmission system, etc) that will lead to a larger percentage increase in the margin.

At present, the concept is at its nascent form and in the future, the analysis will be performed by considering the impact of phase unbalance, tap changing & shunt switching. Furthermore, a theoretical understanding of the behavior of the TDDI as the loading increases is necessary to determine the sufficient conditions to calculate the TDDI at the critical load by just using data at an operating condition. One way forward is to utilize a sensitivity method to estimate the equivalent circuit parameters directly from phasor measurements topology information [18]. The sensitivity based method has been shown to be robust to noise which is a key concern when dealing with measurement based methods. Utilizing the extent of contribution from each sub-system to determine control actions among transmission and distribution networks for effective improvement of the margin is another important analytical step to showcase the methods utility. Testing this method on large transmission systems with several loads replaced by large distribution feeders will ensure that the proposed method can be used for practical systems. Also, more studies need to be conducted on the placement of microPMUs in the distribution system to be able to accurately distinguish between transmission limited and distribution limited systems.

## APPENDIX

The distribution feeder impedances in per unit with voltage base of IEEE 9 bus system are as follows:

| Feeder Config. | Line D0-D1 | Line D1-D2 | Line D1-D3 |
| --- | --- | --- | --- |
| FC1 | 0.33+0.78j | 0.25+0.59j | 0.41+0.98j |
| FC2 | 0.132+1.95j | 0.10+0.089j | 0.164+0.294j |

A load of 4.5 MW and 1.5 MVAR is on Bus D2 and Bus D3. 10 identical feeder configurations are attached in parallel to Bus 5 (D0 corresponds to bus 5), replacing the load of 90 MW and 30 MVAR.

## REFERENCES

[1] M. Glavic and T. Van Cutsem, "A short survey of methods for voltage instability detection," *2011 IEEE Power and Energy Society General Meeting*, San Diego, CA, 2011, pp. 1-8.

[2] K. Vu, M. Begovic, D. Novosel, and M. Saha, "Use of local measurements to estimate voltage-stability margin," *Power Systems, IEEE Transactions on*, vol. 14, no. 3, pp. 1029–1035, Aug 1999.

[3] I. Smon, G. Verbic and F. Gubina, "Local voltage-stability index using tellegen's Theorem," in *IEEE Transactions on Power Systems*, vol. 21, no. 3, pp. 1267-1275, Aug. 2006.

[4] S. Corsi and G. Taranto, "A real-time voltage instability identification algorithm based on local phasor measurements," *Power Systems, IEEE Transactions on*, vol. 23, no. 3, pp. 1271–1279, Aug 2008.

[5] M. Glavic *et al.*, "See It Fast to Keep Calm: Real-Time Voltage Control Under Stressed Conditions," in *IEEE Power and Energy Magazine*, vol. 10, no. 4, pp. 43-55, July-Aug. 2012.

[6] F. Hu, K. Sun, A. Del Rosso, E. Farantatos and N. Bhatt, "Measurement-Based Real-Time Voltage Stability Monitoring for Load Areas," in IEEE Transactions on Power Systems, vol. 31, no. 4, pp. 2787-2798, July 2016.

[7] *Voltage Stability of Power Systems: Concepts, Analytical Tools, and Industry Experience*. IEEE, 1990.

[8] R. B. Prada and L. J. Souza, "Voltage stability and thermal limit: constraints on the maximum loading of electrical energy distribution feeders," *Transm. Distrib. IEE Proc. - Gener.*, vol. 145, no. 5, pp. 573–577, Sep. 1998.

[9] H. Sun, Q. Guo, B. Zhang, Y. Guo, Z. Li, and J. Wang, "Master Slave-Splitting Based Distributed Global Power Flow Method for Integrated Transmission and Distribution Analysis," *IEEE Trans. Smart Grid*, vol. 6, no. 3, pp. 1484–1492, May 2015.

[10] Z. Li, Q. Guo, H. Sun, J. Wang, Y. Xu and M. Fan, "A Distributed Transmission-Distribution-Coupled Static Voltage Stability Assessment Method Considering Distributed Generation," Available for early access online in IEEE Transactions on Power Systems.

[11] R. S. A. Abri, E. F. El-Saadany, and Y. M. Atwa, "Optimal Placement and Sizing Method to Improve the Voltage Stability Margin in a Distribution System Using Distributed Generation," *IEEE Trans. Power Syst.*, vol. 28, no. 1, pp. 326–334, Feb. 2013.

[12] A. Singhal and V. Ajjarapu, "Long-Term Voltage Stability Assessment of an Integrated Transmission Distribution System," *2017 North American Power Symposium (NAPS),* 2017, pp. 1-6. Available online on arXiv: https://128.84.21.199/abs/1711.01374

[13] V. Ajjarapu and C. Christy, "The continuation power flow: a tool for steady state voltage stability analysis," *IEEE Trans. Power Syst.*, vol. 7, no. 1, pp. 416–423, Feb. 1992

[14] R. Bravo and S. Robles, "2012 FIDVR Events Analysis on Valley Distribution Circuits", Report by DER Laboratory at Southern California Edison, July 15, 2013.

[15] M. Jamei, A. Scaglione, C. Roberts, A. McEachern, E. Stewart, S. Peisert, C. McParland, "Online Thevenin Parameter Tracking Using Synchrophasor Data", *2017 IEEE Power and Energy Society General Meeting*, Chicago, IL, 2017, pp. 1-8.

[16] "IEEE 123 Bus Test System," *IEEE Distribution Test Feeders*. [Online]. Available: https://ewh.ieee.org/soc/pes/dsacom/testfeeders/

[17] R. D. Zimmerman, C. E. Murillo-Sanchez, and R. J. Thomas, "MATPOWER: Steady-State Operations, Planning and Analysis Tools for Power Systems Research and Education," Power Systems, IEEE Transactions on, vol. 26, no. 1, pp. 12–19, Feb. 2011.

[18] A. R. R. Matavalam and V. Ajjarapu, "Sensitivity based Thevenin Index with Systematic Inclusion of Reactive Limits", Available for early access online in IEEE Transactions on Power Systems.